\documentclass[a4paper,11pt]{amsart}

\usepackage{amstext}
\usepackage{amsmath}
\usepackage{ifthen}
\usepackage{multirow}
\usepackage{eucal}
\usepackage{amscd}
\usepackage{amssymb}
\usepackage[english]{babel}
\usepackage[T1]{fontenc}
\usepackage[latin1]{inputenc}
\usepackage[all]{xy}
\usepackage{graphicx}
\usepackage{epsfig}

\newenvironment{demo}[1][]{\ifthenelse{\equal{#1}{}}{\noindent \textbf{Proof :\\ \indent}}{\noindent \textbf{Proof #1 :\\ \indent}}}{$\square$\\}

\newtheoremstyle{break}% name
  {}%      Space above, empty = `usual value'
  {}%      Space below
  {\itshape}% Body font
  {}%         Indent amount (empty = no indent, \parindent = para indent)
  {\bfseries}% Thm head font
  {.}%        Punctuation after thm head
  {\newline}% Space after thm head: \newline = linebreak
  {}%         Thm head spec
  
\newtheoremstyle{remarque}
  {}
  {}
  {\slshape}
  {}
  {\bfseries}
  {.}
  {\newline}
  {}

\newtheoremstyle{defin}
  {}
  {}
  {\upshape}
  {}
  {\bfseries}
  {.}
  {\newline}
  {}

\theoremstyle{break}
\newtheorem{theo}{Theorem}[section]
\newtheorem{cor}{Corollary}[section]
\newtheorem{lem}{Lemma}[section]
\newtheorem{prop}{Proposition}[section]
\newtheorem{defi}{Definition}[section]
\theoremstyle{remarque}
\newtheorem{rem}{Remark}[section]
\theoremstyle{defin}
\newtheorem{ex}{Example}[section]

\newcommand{\CC}{\mathbb C}
\newcommand{\QQ}{\mathbb Q}
\newcommand{\ZZ}{\mathbb Z}

\newcommand{\PP}{\mathbb P}

\newcommand{\To}{\longrightarrow}
\newcommand{\abs}[1]{\left\vert#1\right\vert}

\newcommand{\set}[1]{\left\{#1\right\}}
\newcommand{\Forall}[2]{\forall \, #1 \in #2, \:}

\newcommand{\cycl}[1]{\mathcal{C}(#1)}
\newcommand{\hg}[2]{\pi_{#1}(#2)}
\newcommand{\dimm}[1]{\mathrm{dim}(#1)}
\newcommand{\fii}{\varphi}

\newcommand{\wtx}{\widetilde{X}}
\newcommand{\wg}{\widetilde{\Gamma}}
\newcommand{\crochet}[1]{\langle\!\langle\, #1 \,\rangle\!\rangle}

\title{$\Gamma$-reduction for smooth orbifolds}
\author{Benoît \textsc{Claudon}}
\address{Institut \'Elie Cartan \\ UMR 7502\newline\indent Nancy Universit\'e, CNRS, INRIA B.P. 239\newline\indent 54506 Vandoeuvre-l\`es-Nancy Cedex (France)}
\email{Benoit.Claudon@iecn.u-nancy.fr}
\begin{document}

\maketitle

\begin{abstract}
The aim of this short note is to show how to construct a rational Remmert reduction(the $\wg$-reduction) for the universal cover of smooth orbifolds $(X/\Delta)$. Doing this, we are led to introduce some singular Kähler metric on $(X/\Delta)$ adapted to the $\QQ$-divisor $\Delta$.
\end{abstract}

\section*{Introduction}

Let $X$ be a compact Kähler manifold and $\wtx$ be its universal cover. Recall first that the birationnal structure of the latter is partially described by the following result, in the direction of Shafarevich's Conjecture \cite{Sh74} :
\begin{theo}[th. 3.5, p. 264 \cite{Ca94}]\label{existence-gred}
There exists a unique meromorphic fibration (\emph{i.e.} surjective with connected fibers)
$$\gamma_{\wtx}:\wtx\dashrightarrow\Gamma(\wtx)$$
which is almost holomorphic\footnote{the indetermination locus of $\gamma_{\wtx}$ does not surject onto $\Gamma(\wtx)$.}, proper and satisfying the following condition: if $Z\subset\wtx$ is a compact irreducible analytic subset of $\wtx$ passes through a very general point $x\in\wtx$, it is contained in the fiber through $x$
$$Z\subset \gamma_{\wtx}^{-1}\left(\gamma_{\wtx}(x)\right).$$
\end{theo}
\begin{defi}
The fibration $\gamma_{\wtx}$ is called the $\wg$-reduction (or Shafarevich map in the terminology of \cite{K93}) of $\wtx$.
\end{defi}

Here we consider a smooth geometric orbifold $(X/\Delta)$ given by a $\QQ$-divisor
$$\Delta=\sum_{j\in J}(1-\frac{1}{m_j})\Delta_j$$
where $m_j\geq2$ are positive integers and $\textrm{Supp}(\Delta)=\cap_{j\in J}\Delta_j$ is of normal crossings. Following the works of Kato and Namba, we can define a suitable notion of ramified cover for $(X/\Delta)$ and, up to slight modifications for $\Delta$, there exists then a universal cover
$$\pi_{\Delta}:\wtx_{\Delta}\To (X/\Delta)$$
attached to every smooth orbifold $(X/\Delta)$. This paper is devoted to prove the following theorem.
\begin{theo}\label{gred-orbifolde}
Let $(X/\Delta)$ be a smooth Kähler orbifold and $\wtx_{\Delta}$ its universal cover. There exists a unique almost holomorphic proper fibration
$$\tilde{\gamma}_{\Delta}:\wtx_{\Delta}\dashrightarrow\Gamma(\wtx_{\Delta})$$
which satisfies the same condition as in the theorem \ref{existence-gred}: every compact irreducible subvariety of $\wtx_{\Delta}$ passing through a very general point $x\in\wtx_{\Delta}$ is contained in the fiber $\tilde{\gamma}_{\Delta}^{-1}\left(\tilde{\gamma}_{\Delta}(x)\right)$.
\end{theo}

In \cite[th. 11.21]{Ca07}, the $\Gamma$-reduction is constructed for smooth orbifolds $(X/\Delta)$ but the fibration is defined on the orbifold itself (and not on the universal cover). Moreover, the singular metrics introduced here (see section 2) seem to be a natural notion for smooth Kähler orbifolds.\\

Before proving theorem \ref{gred-orbifolde}, we shall start with a brief account of the works of Kato and Namba (ramified covers and fundamental groups for orbifolds). We then introduce some singular Kähler metric adapted to the additionnal orbifold structure : pulling it back, it induces a uniform Kähler structure on the universal cover $\wtx_{\Delta}$ which is sufficient to construct the fibration $\tilde{\gamma}_{\Delta}$.\\

\section{Universal cover for smooth orbifolds}

\subsection{Orbifold fundamental group}

Let us first recall what smoothness means for a pair $(X/\Delta)$.
\begin{defi}\label{defi lisse}
A geometric orbifold $(X/\Delta)$ is said to be smooth if the underlying variety $X$ is a smooth manifold and if the $\QQ$-divisor $\Delta$ has only normal crossings. If in a coordinate patch, the support of $\Delta$ can be defined by an equation
$$\prod_{j=1}^r z_j=0,$$
we will say that these coordinates are adapted to $\Delta$.
\end{defi}

In the category of (smooth) orbifolds, there is a good notion of fundamental group. It is defined in the following way : if $\Delta=\sum_{j\in J}(1-\frac{1}{m_j})\Delta_j$, choose a small loop $\gamma_j$ around each component $\Delta_j$ of the support of $\Delta$ (for instance the boundary of a small disc centered on a smooth point of $\Delta_j$ and transverse to it). Consider now the fundamental group of $X^*=X\backslash \textrm{Supp}(\Delta)$ and its normal subgroup\footnote{it does not depend on the choice of a base point.} generated by the loops $\gamma_j^{m_j}$:
$$\crochet{\gamma_j^{m_j},\,j\in J}\leq\hg{1}{X^*}.$$
\begin{defi}\label{defi pi1 orbifolde}
The fundamental group of $(X/\Delta)$ is defined to be:
$$\hg{1}{X/\Delta}:=\hg{1}{X^*}/\crochet{\gamma_j^{m_j},\,j\in J}.$$
\end{defi}
\begin{ex}\label{exemple courbes}
To illustrate the definition above, consider different orbifold structures on $\PP^1$:
\begin{enumerate}
\item if $\Delta$ has just one point in its support, then $\hg{1}{\PP^1/\Delta}=\set{1}$.
\item for $\Delta=(1-1/m)\set{0}+(1-1/n)\set{\infty}$, we get $\hg{1}{\PP^1/\Delta}=\ZZ/d\ZZ$ where $d=\textrm{gcd}(m,n)$.
\item the double cover of $E\To \PP^1$ ($E$ being an elliptic curve) branches over four points $p_1,\,p_2,\,p_3$ and $p_4$: $\hg{1}{\PP^1/\Delta}$ is thus an extension of $\ZZ^2$ by $\ZZ/2\ZZ$ with $\Delta=\sum_j \frac{1}{2}p_j$.
\end{enumerate}
\end{ex}

\subsection{Branched coverings}

Associated to the fundamental group, there is a notion of ramified cover adapted to an orbifold structure. This is one of the fundamental results of \cite{Ka87} (see also \cite{Na87}).
\begin{defi}\label{defi revetement}
A covering \textbf{branched at most at} \mathversion{bold}$\Delta$\mathversion{normal} is a holomorphic map $\pi:Y\To X$ with
\begin{enumerate}
\item $Y$ normal and connected, $\pi$ with dicrete fibers,
\item $\pi$ induces an étale cover over $X^*$,
\item over $\Delta_j$, the ramification index of $\pi$ is $n_j$, a divisor of $m_j$ (for all $j$),
\item every point $x\in X$ admits some connected neighbourhood $V$ satisfying: every connected component $U$ of $\pi^{-1}(V)$ meets the fiber $\pi^{-1}(x)$ in only one point and the restriction $\pi_{\vert U}:U\To V$ is a proper (finite) map.
\end{enumerate}
We shall say that $\pi$ \textbf{branches at} \mathversion{bold}$\Delta$\mathversion{normal} \textbf{exactly} if $n_j=m_j$ for all $j\in J$.
\end{defi}
As in the absolute case (\emph{i.e.} where $\Delta=\emptyset$), the exists a Galois correspondence between subgroups of $\hg{1}{X/\Delta}$ and coverings of $X$ branched at most at $\Delta$.
\begin{theo}[\cite{Ka87},\cite{Na87}]\label{revetement sous groupe}
If $(X/\Delta)$ is a smooth orbifold, there exists a natural one-to-one correspondence between subgroups $G$ of $\hg{1}{X/\Delta}$ and coverings $\pi:Y\To X$ branched at most at $\Delta$. If the subgroup $G$ is normal (resp. of finite index), the corresponding covering is Galois (resp. finite).
\end{theo}
\begin{rem}
The smoothness assumption is not essential but, in this situation, the local (orbifold) fundamental goups are finite. This finiteness condition is actually the needed one to achieve finiteness as in definition \ref{defi revetement} (4) above.
\end{rem}
The correpondence in theorem \ref{revetement sous groupe} goes in the following way. If $\pi:Y\To X$ is a branched covering (branching at most at $\Delta$), consider the subgroup obtained as the image of the composite morphism:
$$\hg{1}{Y^*}\stackrel{\pi_*}{\To}\hg{1}{X^*}\stackrel{\fii}{\To}\hg{1}{X/\Delta},$$
where $Y^*=\pi^{-1}(X^*)$ and $\fii$ is the natural projection. In the other way, choose a subgroup $G\leq \hg{1}{X/\Delta}$ and consider $G'=\fii^{-1}(G)\leq \hg{1}{X^*}$: it corresponds to an étale covering $\pi:Y^*\To X^*$ ; the finiteness of the local fundamental groups can then be used to complete the covering over the support of $\Delta$. The fact that $\crochet{\gamma_j^{m_j},\,j\in J}\subset G'$  is then equivalent to the ramification condition in the definition \ref{defi revetement} (and $\pi:Y\To X$ is a branched covering according to this definition).
\begin{cor}\label{revetement universel orbifolde}
Corresponding to the trivial subgroup $\set{1}\subset \hg{1}{X/\Delta}$, there exists a simply connected normal complex space $\wtx_{\Delta}$ and a branched covering
$$\pi_{\Delta}:\wtx_{\Delta}\To (X/\Delta)$$
branched at most at $\Delta$. It is called the universal covering of $(X/\Delta)$.
\end{cor}
\begin{rem}\label{singularite}
From the construction itself, it can be easily shown that $\wtx_{\Delta}$ has only quotient singularities (located over the singular locus of $\Delta$
$$\textrm{Sing}(\Delta)=\bigcup_{i\neq j}\Delta_i\cap\Delta_j).$$
It is then a $V$-manifold in the sense of \cite{Sa56}\footnote{for an obvious reason, we will not use the terminology "orbifold" in the sense of Thurston.}.
\end{rem}
\begin{ex}\label{exemple revetement singulier}
Let us consider the following smooth orbifold surface $(\PP^2/\Delta)$ with
$$\Delta=\frac{1}{2}\Delta_1+\frac{1}{2}\Delta_2$$
($\Delta_1$ and $\Delta_2$ two distinct lines in $\PP^2$ meeting in one point $p$). An easy computation shows that
$$\hg{1}{\PP^2/\Delta}=\ZZ/2\ZZ$$
and the universal cover of $(\PP^2/\Delta)$ is thus a double cover $S\To (\PP^2/\Delta)$. Over the point $p$, $S$ has a conic singularity locally given by the quotient
$$\CC^2/\langle(x,y)\mapsto(-x,-y)\rangle.$$
The surface $S$ is actually the cone over the normal rational curve in $\PP^2$.
\end{ex}

\subsection{Regularization of orbifold structures}

To avoid the fact that the universal cover does not necessarily branch exactly at $\Delta$, \cite[§ 11]{Ca07} introduced the notion of regular divisor.
\begin{defi}\label{defi regularisation}
Let $(X/\Delta)$ be a smooth orbifold. Let us denote $d_j$ the order of $\gamma_j$ in the quotient $\hg{1}{X/\Delta}$ ($d_j$ divides $m_j$) and define
$$\Delta_{reg}=\sum_{j\in J}(1-\frac{1}{d_j})\Delta_j.$$
It is called the \textbf{regularization} of $\Delta$. If $\Delta_{reg}=\Delta$, the divisor is said to be \textbf{regular}.
\end{defi}
\begin{ex}\label{exemple regularisation}
In the example \ref{exemple courbes}, the regularization of $\Delta$ is given by
\begin{enumerate}
\item $\Delta_{reg}=\emptyset$.
\item $\Delta_{reg}=(1-1/d)\set{0}+(1-1/d)\set{\infty}$ with $d=\textrm{gcd}(m,n)$ (in particular, $\Delta_{reg}=\emptyset$ if $d=1$).
\item $\Delta_{reg}=\Delta$.
\end{enumerate}
\end{ex}
By construction, we can see that
$$\hg{1}{X/\Delta_{reg}}=\hg{1}{X/\Delta}.$$
Moreover, the integers $d_j$ are also the ramification indices of the branched covering $\pi_{\Delta}:\wtx_{\Delta}\To (X/\Delta)$.
\begin{prop}\label{remplacement regularise}
If $(X/\Delta)$ is a smooth orbifold, the orbifold $(X/\Delta_{reg})$ has the same fundamental group and the same universal cover and the covering map 
$$\pi_{\Delta}:\wtx_{\Delta}\To (X/\Delta_{reg})$$
branches exactly at $\Delta_{reg}$.
\end{prop}
Up to regularization (which does not change fundamental group and universal cover), we can then assume that $\Delta$ is regular and that the universal covering branches exactly at $\Delta$.\\

\noindent\textbf{Assumption :} In the rest of the paper, we assume the $\QQ$-divisor $\Delta$ to be \textbf{regular}.

\section{$\wg$-reduction of $\wtx_{\Delta}$}

\subsection{Singular Kähler metrics}

To apply the construction of \cite{Ca94}, we only need a sufficiently uniform Kähler metric on the universel cover $\wtx_{\Delta}$ of a smooth Kähler\footnote{an orbifold is said to be Kähler if the underlying manifold is so.} orbifold $(X/\Delta)$. If $\pi_{\Delta}$ is the covering map et $\omega$ any Kähler metric on $X$, $\pi_{\Delta}^*\omega$ is a closed non negative $(1,1)$-form on $\wtx_{\Delta}$. Unfortunately, over $\Delta$, $\pi_{\Delta}$ is not a local isomorphism and $\pi_{\Delta}^*\omega$ is degenerate. We have to introduce some singularity (concentrated on $\Delta$) to balance the ramification of $\pi_{\Delta}$.
\begin{prop}\label{metrique orbifolde}
Let $(X/\Delta)$ be a smooth Kähler orbifold with $\Delta=\sum_{j\in J}(1-1/m_j)\Delta_j$. Let $\omega$ be any Kähler metric on $X$, let $C>0$ be a real number and $s_j\in H^0(X,\mathcal{O}_X(\Delta_j))$ be a section defining $\Delta_j$. Consider the following expression:
$$\omega_{\Delta}=C\omega+\sum_{j\in J}i\partial\overline{\partial}\abs{s_j}^{2/m_j}$$
where $\abs{\cdot}$ is any smooth metric on $\mathcal{O}_X(\Delta_j)$ (for each $j$). If $C$ is large enough, the above formula defines a closed positive $(1,1)$-current (smooth away from $\Delta$) satisfying moreover
$$\omega_{\Delta}\geq\omega$$
in the sense of currents.
\end{prop}
\begin{rem}\label{modele local}
A look at the local model enlights the previous proposition. Consider $\CC^n$ with the orbifold divisor given by the equation
$$\prod_{j=1}^n z_j^{1-1/m_j}=0$$
(with eventually $m_j=1$ for some $j$). The sections $s_j$ are simply the coordinates $z_j$ and a simple computation gives
\begin{align*}
\omega_{\Delta}&=\omega_{eucl}+\sum_{j=1}^n i\partial\overline{\partial}\abs{z_j}^{2/m_j}\\
&=\omega_{eucl}+\sum_{j=1}^n\frac{idz_j\wedge d\overline{z_j}}{m_j^2\abs{z_j}^{2(1-1/m_j)}}
\end{align*}
Moreover, the uniformization is given by
$$\pi:\left\{\begin{array}{ccc}\CC^n & \To & \CC^n\\
(t_1,\dots,t_n) & \mapsto & (z_1,\dots,z_n)=(t_1^{m_1},\dots,t_r^{m_r},t_{r+1},\dots,t_n)
\end{array}\right.$$
and, in this chart, the above expression becomes
\begin{align*}
\pi^*(\omega_{\Delta})&=\sum_{j=1}^n\frac{id(t_j^{m_j})\wedge d(\overline{t_j}^{m_j})}{\abs{t_j^{m_j}}^{2(1-1/m_j)}}+\pi^*(\omega_{eucl})\\
&=\sum_{j=1}^n\frac{im_j^2t_j^{m_j-1}dt_j\wedge \overline{t_j}^{m_j-1}d\overline{t_j}}{\abs{t_j}^{2(m_j-1)}}+\pi^*(\omega_{eucl})\\
&=\sum_{j=1}^n im_j^2(1+\abs{t_j}^{2(m_j-1)})dt_j\wedge d\overline{t_j}.
\end{align*}
In the uniformization $\pi$, the $(1,1)$-form $\omega_{\Delta}$ becomes a genuine Käher metric.
\end{rem}
\begin{demo}[of proposition \ref{metrique orbifolde}]
We only need to check that each $\abs{s_j}^{2/m_j}$ is a quasi-psh function on $X$. In adapted coordinates, the sections $s_j$ are given by $z_j$ but we have to take care of the weights of the metrics (on $\mathcal{O}_X(\Delta_j)$):
$$\abs{s_j}^{2/m_j}=f_j\abs{z_j}^{2/m_j}$$
with $f_j$ a smooth positive function. We then get
\begin{align}\label{calcul metrique}
i\partial\overline{\partial}\abs{s_j}^{2/m_j}&=if_j\partial\overline{\partial}\abs{z_j}^{2/m_j}+ i\partial f_j\wedge\overline{\partial}\abs{z_j}^{2/m_j}\nonumber\\
&+i\partial\abs{z_j}^{2/m_j}\wedge\overline{\partial}f_j+ i\abs{z_j}^{2/m_j}\partial\overline{\partial}f_j.
\end{align}
The following identity
\begin{align*}
0\leq i\partial(\abs{z_j}^{2/m_j}+f_j)\wedge&\overline{\partial}(\abs{z_j}^{2/m_j}+f_j)=i\partial f_j\wedge\overline{\partial}f_j+i\partial f_j\wedge\overline{\partial}\abs{z_j}^{2/m_j}\\
&+i\partial\abs{z_j}^{2/m_j}\wedge\overline{\partial}f_j+i\partial \abs{z_j}^{2/m_j}\wedge\overline{\partial}\abs{z_j}^{2/m_j},
\end{align*}
gives the inequality (in the sense of currents):
\begin{align}\label{minoration metrique1}
i\partial f_j\wedge\overline{\partial}\abs{z_j}^{2/m_j}+i\partial\abs{z_j}^{2/m_j}\wedge\overline{\partial}f_j 
&\geq-i\partial f_j\wedge\overline{\partial}f_j-i\partial \abs{z_j}^{2/m_j}\wedge\overline{\partial}\abs{z_j}^{2/m_j}\nonumber\\
&\geq-i\partial f_j\wedge\overline{\partial}f_j-\frac{i\abs{z_j}^{2/m_j}dz_j\wedge d\overline{z_j}}{\abs{z_j}^{2(1-1/m_j)}}.
\end{align}
Since $f_j$ is smooth, there exists a constant $C_j>0$ such that (locally)
\begin{equation}\label{minoration metrique2}
i\abs{z_j}^{2/m_j}\partial\overline{\partial}f_j-i\partial f_j\wedge\overline{\partial}f_j\geq -C_j\omega
\end{equation}
Combining (\ref{calcul metrique}), (\ref{minoration metrique1}) and (\ref{minoration metrique2}) gives
\begin{align*}
\omega_{\Delta}&\geq (C-\sum_{j}C_j)\omega+\sum_{j}\frac{(f_j-\abs{z_j}^{2/m_j})idz_j\wedge d\overline{z_j}}{\abs{z_j}^{2(1-1/m_j)}}\\
&\geq(C-\sum_{j}C_j)\omega\,\textrm{ on a neighbourhood of }0.
\end{align*}
Since $X$ is compact, it is covered by a finite number of such small balls and we can choose $C$ large enough to achieve positivity for $\omega_{\Delta}$.
\end{demo}

\noindent Pulling back this singular metric to the universal cover, we get the needed uniform metric.
\begin{prop}
Choose $\omega_{\Delta}$ as in the proposition \ref{metrique orbifolde}. The pull-back
$$\widetilde{\omega_{\Delta}}=\pi_{\Delta}^*(\omega_{\Delta})$$
defines a Kähler metric on $\wtx_{\Delta}$ as a $V$-manifold. This means that the (1,1)-form $\widetilde{\omega_{\Delta}}$ is continuous on $\wtx_{\Delta}$, $C^{\infty}$ on the non-singular locus of $\wtx_{\Delta}$ and its lift to any local uniformization extends to a smooth invariant metric\footnote{a $V$-manifold is locally a quotient $U/G$ where $U\subset\CC^n$ is an open subset and $G$ a finite group acting on $U$ ; in this setting, a lift is just the pull-back to $U$.}.
\end{prop}
\begin{demo}
The local computation made in the remark \ref{modele local} shows exactly that the singular part of the metric $\omega_{\Delta}$  balances the ramification of $\pi_{\Delta}$.
\end{demo}

\subsection{Proof of theorem \ref{gred-orbifolde}}

We can now prove the existence of the $\wg$-reduction for $\wtx_{\Delta}$. We will need some compactness properties ot the cycles space $\cycl{\wtx_{\Delta}}$ constructed in \cite{Ba75}. To construct the $\wg$-reduction, we will apply the fundamental result of \cite{Ca81}. For the convenience of the reader, we restate it here.
\begin{theo}[th. 1, p. 189 \cite{Ca81}]\label{T quotient}
Let $Y$ be a normal analytic space, $T$ an irreducible component of $\cycl{Y}$,
$$G_T=\set{(y,Z_t)\in Y\times T\vert \, y\in\textrm{Supp}(Z_t)}$$
the incidence graph of the family of cycles parametrized by $T$ and let us denote by
$$q:G_T\To Y\textrm{ et }r:G_T\To T$$
the corresponding projections. If $q$ is surjective (\emph{i.e.} the cycles $(Z_t)_{t\in T}$ cover $Y$) and proper and if the generic fiber of $r$ is irreducible (\emph{i.e.} $Z_t$ is irreducible for $t\in T$ generic), there exists an almost holomorphic proper fibration
$$g_T:Y\dashrightarrow Q_T$$
whose fiber through a generic point $y\in Y$ is the equivalence class generated by this family of cycles (two points are said to be equivalent if there is a connected (finite) chain of cycles of the family containing them).
\end{theo}
\noindent We thus only need to prove the following
\begin{prop}\label{proprete cycle}
If $T\subset\cycl{\wtx_{\Delta}}$ is an irreducible component and $q:G_T\To \wtx_{\Delta}$ (and $r:G_T\To T$) is the corresponding projection, then $q$ is proper.
\end{prop}
\noindent This proposition is a straightforward consequence of the following lemma which restates the fact that the geometry associated to the Kähler form $\widetilde{\omega_{\Delta}}$ is uniform (see \cite{B92} for the metric structure of $V$-manifold).
\begin{lem}\label{lemme Bishop}
There exist some constants $r,\delta>0$ such that: for every irreducible compact subvariety $Z\subset\wtx_{\Delta}$ and every $z\in Z$,
$$\mathrm{Vol}_{\widetilde{\omega_{\Delta}}}(Z\cap B(z,r)):=\int_{Z\cap B(z,r)} \widetilde{\omega_{\Delta}}^{\dimm{Z}}\geq \delta$$
where $B(z,r)$ is the ball of radius $r$ centered at $z$ for the distance $d_{\Delta}$ induced by $\widetilde{\omega_{\Delta}}$ on $\wtx_{\Delta}$.
\end{lem}
\noindent Here, the volume of an irreducible subvariety is computed as follows: if $Z$ is contained in a chart $p:U\To V\simeq U/G$, $$\int_Z\widetilde{\omega_{\Delta}}^{\dimm{Z}}=\frac{1}{\abs{G}}\int_{p^{-1}(Z)}p^*(\widetilde{\omega_{\Delta}})^{\dimm{Z}}$$
as usually (note that the metric $p^*(\widetilde{\omega_{\Delta}})$ is now smooth) and we use a partition of unity to deal with the general case.\\

\begin{demo}[of the proposition \ref{proprete cycle}]
Let $T$ be an irreducible component of $\cycl{\wtx_{\Delta}}$ and $K$ a compact of $\wtx_{\Delta}$. Since $\widetilde{\omega_{\Delta}}$ is a Kähler metric, the volume of the cycles parametrized by $T$ is constant; let us denote it by $v$. Consider the following compact subset of $\wtx_{\Delta}$:
$$\hat{K}=\set{x\in \wtx_{\Delta}\vert\, d_{\Delta}(x,K)\leq M}\textrm{ with }M>r\left\lceil\frac{v}{\delta}\right\rceil.$$
The lemma \ref{lemme Bishop} can then be rephrased in the following way:
$$q^{-1}(K)\subset r^{-1}(r(q^{-1}(\hat{K}))).$$
But Bishop's theorem asserts the compactness of $r(q^{-1}(\hat{K}))$ (see \cite{Li78}) and the projection $r$ is always proper. From this we deduce the compactness of $q^{-1}(K)$ (and the properness of $q$).
\end{demo}
\newpage
\noindent We can now finish the\\

\begin{demo}[of the theorem \ref{gred-orbifolde}]
Let us $d$ denote the smallest integer such that there exists an almost holomorphich proper fibration
$$f:\wtx_{\Delta}\To V$$
with $\dimm{V}=d$. Assume that the maximality property of the fibers of $f$ (stated in theorem \ref{gred-orbifolde}) is not satisfied: for a very general point $x\in\wtx_{\Delta}$, there exists a compact irreducible subvariety not contained in the fiber through $x$. Since $\cycl{\wtx_{\Delta}}$ is a second countable space, there exists an irreducible component $T$ of $\cycl{\wtx_{\Delta}}$ such that the family $(U_t)_{t\in T}$ of cycles parametrized by $T$ satisfy
$$\Forall{t}{T}\dimm{f(U_t)}>0$$
and $U_t$ is irreducible for a generic $t\in T$. Thanks to proposition \ref{proprete cycle}, we can apply theorem \ref{T quotient} to the family of \emph{compact} cycles
$$Z_t=f^{-1}(f(U_t))$$
parametrized (meromorphically) by $T$. The corresponding quotient
$$g_T:\wtx_{\Delta}\dashrightarrow Q_T$$
is an almost holomorphic proper fibration whose fibers are strictly contained in the one of $f$; from this we deduce
$$\dimm{Q_T}<d=\dimm{V},$$
which contradicts the minimality of $d$.\\
The uniqueness of the $\wg$-reduction follows in the same way: if $f$ and $g$ are two such fibrations and $x$ a sufficiently general point, the fiber of $f$ through $x$ is contained in the fiber of $g$ (maximality of the fibers of $g$). Reversing the order of $f$ and $g$, we get the other inclusion and the two fibrations have the same fibers.
\end{demo}
\begin{rem}
The preceeding construction can be made with an arbitrary non compact branched covering $Y\To (X/\Delta)$ (branched at most at $\Delta$). The regularization $\Delta_Y$ must be adapted to the previous map and it yields
$$\Delta\geq \Delta_Y \geq \Delta_{reg}$$
(the fundamental group is thus unchanged). The pull-back of the singular metric $\omega_{\Delta_Y}$ to $Y$ gives the right object to apply the construction described above.
\end{rem}

To conclude, we would like to recall the results of \cite[th. 11.21]{Ca07} and to compare both fibrations.
\begin{theo}\label{gred-fred}
Let $(X/\Delta)$ be a smooth Kähler orbifold. There exists a unique almost holomorphic fibration
$$\gamma_{(X/\Delta)}:(X/\Delta)\dashrightarrow\Gamma(X/\Delta)$$
satisfying both following properties:
\begin{enumerate}
\item $\hg{1}{X_y/\Delta_{X_y}}_{(X/\Delta)}:=\textrm{Im}\left(\hg{1}{X_y/\Delta_{X_y}}\To \hg{1}{X/\Delta}\right)$ is finite for $y\in \Gamma(X/\Delta)$ general,
\item if $g:(V/\Delta_V)\To (X/\Delta)$ is a divisible orbifold morphism from a smooth compact orbifold $(V/\Delta_V)$ such that $g(V)$ meets $X_y$ (for $y\in \Gamma(X/\Delta)$ generic) and if
$$\textrm{Im}\left(\hg{1}{V/\Delta_V}\stackrel{g_*}{\To} \hg{1}{X/\Delta}\right)$$
is finite, then $g(V)$ is contained in $X_y$.
\end{enumerate}
\end{theo}
\noindent Recall that a divisible orbifold morphism
$$g:(V/\Delta_V)\To (X/\Delta)$$
induces a well-defined morphism at the level of fundamental groups:
$$g_*:\hg{1}{V/\Delta_V}\To\hg{1}{X/\Delta}.$$
See \cite[§ 2.2 and § 11.1]{Ca07} for the notions involved.\\

In the tame situation of a smooth subvariety meeting $\Delta$ transversally\footnote{by this we mean that $(V/i^*\Delta)$ should be a smooth orbifold.} (which is the case of the general fiber of a fibration), the inclusion
$$i:(V/\Delta_{\vert V})\hookrightarrow (X/\Delta)$$
is clearly a divisible orbifold morphism. The following lemma shows us that the fibers of $\tilde{\gamma}_{\Delta}$ are the connected components of the inverse images by the covering map $\pi_{\Delta}$ of the fibers of $\gamma_{(X/\Delta)}$.
\begin{lem}
Let $V$ be a smooth subvariety of $X$ meeting $\Delta$ transversally and $\pi_{\Delta}:\wtx_{\Delta}\To (X/\Delta)$ the universal cover of $(X/\Delta)$. Both conditions are equivalent :
\begin{enumerate}
\item $\hg{1}{V/\Delta_V}_{(X/\Delta)}$ is a finite group.
\item each connected component of $\pi_{\Delta}^{-1}(V)$ is compact.
\end{enumerate}
\end{lem}
\begin{demo}
Let $Z$ be such a connected component. Restricting $\pi_{\Delta}$ to $Z$ yields a branched covering
$$p:Z\To (V/\Delta_V)$$
wich corresponds to the subgroup
$$G=\textrm{Ker}\left(\hg{1}{V/\Delta_V}\To\hg{1}{X/\Delta}\right).$$
According to theorem \ref{revetement sous groupe}, $p$ is finite if and only if $G$ is of finite index in $\hg{1}{V/\Delta_V}$. Equivalently, $Z$ is compact if and only if $\hg{1}{V/\Delta_V}_{(X/\Delta)}$ is finite.
\end{demo}

\noindent\textbf{Aknowledgements.} I am very grateful to Frédéric Campana for introducing me to this subject. For this and also for many interesting discussions, I would like to thank him. 

\newpage

%Bibliographie---------------------------------------------------------------
\bibliographystyle{amsalpha}
\bibliography{myref}

\providecommand{\bysame}{\leavevmode\hbox to3em{\hrulefill}\thinspace}
\providecommand{\MR}{\relax\ifhmode\unskip\space\fi MR }
\providecommand{\MRhref}[2]{%
  \href{http://www.ams.org/mathscinet-getitem?mr=#1}{#2}
}
\providecommand{\href}[2]{#2}
\begin{thebibliography}{Nam87}

\bibitem[Bar75]{Ba75}
D.~Barlet, \emph{Espace analytique r\'eduit des cycles analytiques complexes
  compacts d'un espace analytique complexe de dimension finie}, Fonctions de
  plusieurs variables complexes, II (S\'em. Fran\c{c}ois Norguet, 1974--1975),
  Springer, Berlin, 1975, pp.~1--158. Lecture Notes in Math., Vol. 482.

\bibitem[Bor92]{B92}
J.~Borzellino, \emph{Riemannian geometry of orbifolds}, Ph.D. thesis,
  University of California, Los Angeles, 1992.

\bibitem[Cam81]{Ca81}
F.~Campana, \emph{Cor\'eduction alg\'ebrique d'un espace analytique faiblement
  k\"ahl\'erien compact}, Invent. Math. \textbf{63} (1981), no.~2, 187--223.

\bibitem[Cam94]{Ca94}
\bysame, \emph{Remarques sur le rev\^etement universel des vari\'et\'es
  k\"ahl\'eriennes compactes}, Bull. Soc. Math. France \textbf{122} (1994),
  no.~2, 255--284.

\bibitem[Cam07]{Ca07}
\bysame, \emph{Orbifoldes sp\'eciales et classifications bim\'eromorphes des
  vari\'et\'es k\"ahl\'eriennes compactes}, preprint arXiv:0705.0737, 2007.

\bibitem[Kat87]{Ka87}
Mitsuyoshi Kato, \emph{On uniformizations of orbifolds}, Homotopy theory and
  related topics (Kyoto, 1984), Adv. Stud. Pure Math., vol.~9, North-Holland,
  Amsterdam, 1987, pp.~149--172.

\bibitem[Kol93]{K93}
J.~Koll{\'a}r, \emph{Shafarevich maps and plurigenera of algebraic varieties},
  Invent. Math. \textbf{113} (1993), no.~1, 177--215.

\bibitem[Lie78]{Li78}
D.~I. Lieberman, \emph{Compactness of the {C}how scheme: applications to
  automorphisms and deformations of {K}\"ahler manifolds}, Fonctions de
  plusieurs variables complexes, III (S\'em. Fran\c cois Norguet, 1975--1977),
  Lecture Notes in Math., vol. 670, Springer, Berlin, 1978, pp.~140--186.

\bibitem[Nam87]{Na87}
M.~Namba, \emph{Branched coverings and algebraic functions}, Pitman Research
  Notes in Mathematics Series, vol. 161, Longman Scientific \& Technical,
  Harlow, 1987.

\bibitem[Sat56]{Sa56}
I.~Satake, \emph{On a generalization of the notion of manifold}, Proc. Nat.
  Acad. Sci. U.S.A. \textbf{42} (1956), 359--363.

\bibitem[Sha74]{Sh74}
I.~R. Shafarevich, \emph{Basic algebraic geometry}, Springer-Verlag, New York,
  1974, Translated from the Russian by K. A. Hirsch, Die Grundlehren der
  mathematischen Wissenschaften, Band 213.

\end{thebibliography}

\end{document}